\documentclass[letterpaper, 10 pt, conference]{ieeeconf}  

\IEEEoverridecommandlockouts                              
\usepackage{graphics}
\usepackage{epsfig}
\usepackage{amsmath}
\usepackage{amssymb}
\usepackage{color}
\usepackage{cite}
\usepackage{pst-all} 
\usepackage{relsize} 
\usepackage{ifpdf}

\def\N{{\mathbb{N}}}
\def\R{{\mathbb{R}}}
\def\Q{{\mathbb{Q}}}

\def\Black{\textcolor{black}}

\definecolor{karl}{rgb}{.9,.1,.4} 
\definecolor{juergen}{rgb}{.1,.4,.9} 
\def\karl{\Black}
\def\juergen{\Black}
\def\kawo{\Black}

\def\state{{x}}

\def\statestar{{\state^\star}}
\def\stateset{\karl{{\mathbb{R}^n}}}
\def\constraintstateset{\karl{X}}
\def\control{{u}}
\def\controlstar{{\u^\star}}
\def\controlset{\karl{{\mathbb{R}^m}}}
\def\constraintcontrolset{\karl{U}}
\def\admissiblecontrolset{\mathcal{U}}
\def\feedback{{\mu}}

\def\stagecost{{\ell}}

\def\updateinstant{{\sigma}}

\def\tsc{{\tau_{sc}}}
\def\tca{{\tau_{ca}}}
\def\tc{{\tau_{c}}}
\def\hmpc{{\hspace*{-0.375mm}}}
\def\hinfty{{\hspace*{-0.375mm}}}
\def\z{{\delta}}
\def\zt{{\tilde{\delta}}}

\def\x{{x}}
\def\X{{\stateset}}
\def\Xb{{\constraintstateset}}
\def\u{{\control}}
\def\U{{\controlset}}
\def\Ub{{\constraintcontrolset}}
\def\Uc{{\admissiblecontrolset}}

\def\xs{{\x^\star}}
\def\xh{{\hat{\x}}}
\def\us{{\u^\star}}
\def\l{{\stagecost}}

\def\tsc{{\tau_{sc}}}
\def\tca{{\tau_{ca}}}
\def\tc{{\tau_{c}}}
\def\hmpc{{\hspace*{-0.375mm}}}
\def\hinfty{{\hspace*{-0.375mm}}}
\def\hv{\vspace*{2.5mm}}
\def\m{{m}}
\def\mt{{\widetilde{m}}}
\def\s{{\tau}}

\newtheorem{theorem}{Theorem}
\newtheorem{definition}[theorem]{Definition}
\newtheorem{proposition}[theorem]{Proposition}

\newtheorem{remark}[theorem]{Remark}
\newtheorem{corollary}[theorem]{Corollary}
\newtheorem{assumption}[theorem]{Assumption}

\title{\LARGE \bf
Ensuring Stability in Networked Systems with Nonlinear MPC for Continuous Time Systems
}

\author{Lars Gr\"{u}ne$^{1}$, J\"{u}rgen Pannek$^{2}$, and Karl Worthmann$^{1}$%
\thanks{$^{1}$L.~Gr\"{u}ne and K.~Worthmann are with the Mathematical Institute, University of Bayreuth, 95440 Bayreuth, Germany
        {\tt\small (lars.gruene, karl.worthmann@uni-bayreuth.de)}.}%
\thanks{$^{2}$J.~Pannek is with the Faculty of Aerospace Engineering, University of the Federal Armed Forces Munich, 85577 Munich/Neubiberg, Germany {\tt\small (juergen.pannek@unibw.de)}.}%
}

\begin{document}

\maketitle
\thispagestyle{empty}
\pagestyle{empty}

\begin{abstract}
	For networked systems, the control law is typically subject to network flaws such as delays and packet dropouts. Hence, the time in between updates of the control law varies unexpectedly. Here, we present a stability theorem for nonlinear model predictive control with varying control horizon in a continuous time setting without stabilizing terminal constraints or costs. It turns out that stability can be concluded under the same conditions as for a (short) fixed control horizon.
\end{abstract}

\section{Introduction}
In recent years, networked control systems (NCS) received growing popularity due to their lower implementation costs and greater interoperability compared to standard control system\karl{s, cf.~\cite{GuptaChow2010}}. 
NCS designs have been implemented in different areas such as robotics, automotive and aeronautical applications, \karl{see, e.g., \cite{LeenHeffernan2002,OdaDoiWakata2001}}. On the backside, however, the stability and performance analysis of feedbacks designed for networked control systems is more complex, cf.~\cite{vanDeWouwNesicHeemels2012}.

\karl{In this paper we investigate} stability and performance of a nonlinear model predictive controller (MPC) in a prediction consistent network scheme. By now, MPC has been understood quite well even in the nonlinear case, see, \karl{e.g. \cite{RawlingsMayne2009, GruenePannek2011}}. The beauty of the method lies in its simplicity: \karl{First, an optimal control is computed over a finite optimization horizon}
. Then, a portion of this control is implemented --- throughout this paper the length of this portion is called control horizon --- and last the optimization horizon is shifted in time rendering the method iteratively applicable. Prediction consistency formalizes the equivalence of the control input history for the actuator and the controller. Due to computing and transmission times, delays and dropouts, this equivalence is subject to different time levels within the NCS components. Hence, the control histories do not coincide automatically but this property can be forced using certain communication schemes, cf.~\cite{Bemporad1998, GruenePannekWorthmann2009b, PolushinLiuLung2008}. At the center of these schemes is the requirement of long control sequences being transmitted from the controller to the actuator. \karl{MPC is ideally suited to obtain such a sequence} since it not only delivers such a control by construction, but its internal model can also be used to compensate for time delays between the NCS components. 

The prediction consistent approach naturally leads to the fact that the control horizon varies over time. The goal of this paper is to provide the theoretical foundation \juergen{for} ensuring stability of MPC without stabilizing terminal costs or constraints for such varying control horizons in a nonlinear continuous time setting. \juergen{To this end, we extend results from \cite{RebleAllgoewer2011} where a stability condition was introduced for a fixed, typically short control horizon to the networked context. In particular, the condition allows us to prove stability for a large range of possibly time varying control horizons \karl{in analogy to \cite{GruenePannekWorthmann2009ECC} where the same assertion was shown for discrete time systems}.} 
\juergen{This main result is complemented by an} analysis of structural properties \karl{with respect to the overshoot bound and the decay rate of the imposed controllability condition}.

The paper is organized as follows. In Section \ref{SectionPreliminaries} the detailed problem formulation is stated. The results from \cite{RebleAllgoewer2011} are briefly summarized in Section \ref{SectionRecap} and then extended to the time varying case in Section \ref{SectionResults}. Finally, we show numerical results to illustrate our results and draw some conclusions.

\section{Problem setting}\label{SectionPreliminaries}

We consider \karl{a} nonlinear \karl{control} systems governed by a differential equation
\begin{equation}\label{NotationSystemDynamics}
	\dot{\x}(t) = f(\x(t),\u(t)).
\end{equation}
Here\karl{,} $\x(t)$ and $\u(t)$ denote the state and the control at time $t \geq 0$, respectively. The system dynamics is given by \karl{$f: \X \times \U \rightarrow \X$, 
s}tate and control constraints are represented by suitable subsets $\Xb \subset \X$ and $\Ub \subset \U$, respectively. The trajectory which emanates from initial state \karl{$\x_{0}$} and is manipulated by the control function $\u:\R_{\geq 0} \rightarrow \U$ is denoted by \karl{$\x(t;\x_{0},u)$.} We call a control function $\u$ admissible for $\x$ on the interval $[0,T)$ if the conditions
\begin{equation}\label{NotationFeasibility}
	\karl{\x(t;x,u)} \in \Xb, t \in [0,T]\karl{,} \quad\text{and}\quad \u(t) \in \Ub, t \in [0,T)\karl{,}
\end{equation}
hold which \karl{is denoted} by $\u \in \Uc_{\x}([0,T))$. Furthermore, a control function $\u:\R_{\geq 0} \rightarrow \U$ is said to be admissible for $\x$ if, for each $T > 0$, $\u \in \Uc_{\x}([0,T))$ holds for its restriction to $[0,T)$. Then, we write $\u \in \Uc_{\x}([0,\infty))$. 

Here, we consider a networked situation, that is System \eqref{NotationSystemDynamics} is connected to an external controller via a network which may be subject to delays and packet dropouts. Within the network, we suppose the clocks at sensor, controller and actuator to be synchronized, see \cite[Section III.C]{VaruttiFindeisen2009} for a relaxation of this assumption. These clocks are denoted by $t_s$, $t_c$ and $t_a$\karl{,} respectively. 

Our goal is to stabilize System \eqref{NotationSystemDynamics} at an equilibrium $\xs$ for which we suppose that a control input $\us$ exists satisfying $f(\xs,\us) = 0$. To solve this task in an optimal fashion, we introduce continuous running costs $\l:\X \times \U \rightarrow \R_{\geq 0}$ which satisfy
\begin{equation*}
	\l(\xs,\us) = 0 \quad\text{and}\quad \inf_{\u \in \U} \l(\x,\u) \geq \underline{\eta}(\| \x - \xs \|), \x \neq \xs\karl{,}
\end{equation*}
for a $\mathcal{K}_\infty$-function $\underline{\eta}$. {As usual, a continuous function $\eta: \R_{\geq} \rightarrow \R_{\geq 0}$ is said to belong to $\mathcal{K}_\infty$ if $\eta(0) = 0$, it is strictly increasing, and unbounded.
Utilizing the running costs $\l$, we define the cost functional
\begin{equation*}
	J_\infty(\x,\u) := \int_0^\infty \l(\karl{\x(t;\x,\u}),\u(t))\, dt
\end{equation*}
\karl{which we wish to minimize for a given initial value $x = \x_0 \in \X$}
. \karl{The corresponding optimal value function is denoted} by
\begin{equation}\label{NotationOptimalValueFunction}
	V_\infty(\x) := \inf_{\u \in \Uc_{\x}([0,\infty))} J_\infty(\x,\u).
\end{equation}
Since optimal control problems on an infinite time horizon are computationally hard, we use model predictive control (MPC) in order to approximately solve \karl{this 
task}. To this end, we firstly set $\xh:=\x_0$. Then, our MPC scheme consists of the following three steps:
\begin{itemize}
	\item Compute a minimizing control function $\u^\star$ for the optimal control problem on a truncated and, thus, finite prediction horizon $T$ depending on $\xh$, i.e.
		\begin{equation}\label{NotationOptimalValueFunctionFiniteHorizon}
			\min_{\u \in \karl{\Uc}_{\xh}[0,T)} \hspace*{-.25mm} J_T(\xh,\u) = \hspace*{-.25mm} \min_{\u \in \karl{\Uc}_{\xh}[0,T)} \int_0^T \hspace*{-.5mm} \l(\x(t;\xh,\u),\u(t)) dt.
		\end{equation}
	\item Define the MPC feedback law $\mu_T: [0,T) \times \X \rightarrow \U$ by $\mu_{T}(t,\xh) := \u^\star(t)$, $t \in [0,T)$. Then, for given control horizon $\delta \in [0,T]$, implement the first piece $\u^\star(t)|_{t \in [0,\delta)}$ of the computed control at the plant in order to obtain
		\begin{equation*}
			\x_{\mu_{T,\delta}}(\delta;\xh) = \x(\delta;\xh,\u^\star).
		\end{equation*}
		Whenever we want to reflect the control horizon $\delta$ in the notation we write $\mu_{T,\delta} = \mu_{T}$.
	\item Shift the prediction (optimization) horizon forward in time by $\delta$ and obtain the new state measurement $\xh$ which coincides in nominal MPC with $\x_{\mu_{T,\delta}}(\delta;\xh)$.
\end{itemize}
Iterative application of this procedure generates a solution on the infinite time horizo\karl{n. 
T}he resulting input signal and trajectory  at time $t$ are denoted by $\mu^{\text{MPC}}(t;x_0)$ and $\x^{\text{MPC}}(t;x_0)$, respectively. Here, we tacitly assume that Problem \eqref{NotationOptimalValueFunctionFiniteHorizon} is solvable and the minimum is attained in each step of the proposed MPC algorithm. For a detailed discussion of feasibility issues we refer to \cite{RawlingsMayne2009,GruenePannek2011}. 

In order to compensate for the mentioned network flaws, we utilize a communication protocol which transmits time stamped information between the network components: For one, the sensor transmits the latest state measurement $\xh$ together with the measurement time instant $t_s$ to the controller. \karl{S}econdly, the controller sends predicted control functions $\feedback_{T, \Delta, \updateinstant} = \controlstar_{[0, \Delta)}$ with $\delta \leq \Delta \leq T$ where $\updateinstant \in \R$ corresponds to the time at which the control should be applied. Here, we assume that the controller computes the input functions $\feedback_{T, \Delta, \updateinstant}$ at a predefined sampling rate (typically with sampling time $<\hspace{-1mm}<\Delta$) based on the most recent measurement available in the buffer. To store these controls, we add a buffer to the actuator and to the controller, cf.~Fig.~\ref{fig:networked system}.

\ifpdf
\begin{figure}[!ht]
	\begin{center}
		\includegraphics[width=0.40\textwidth]{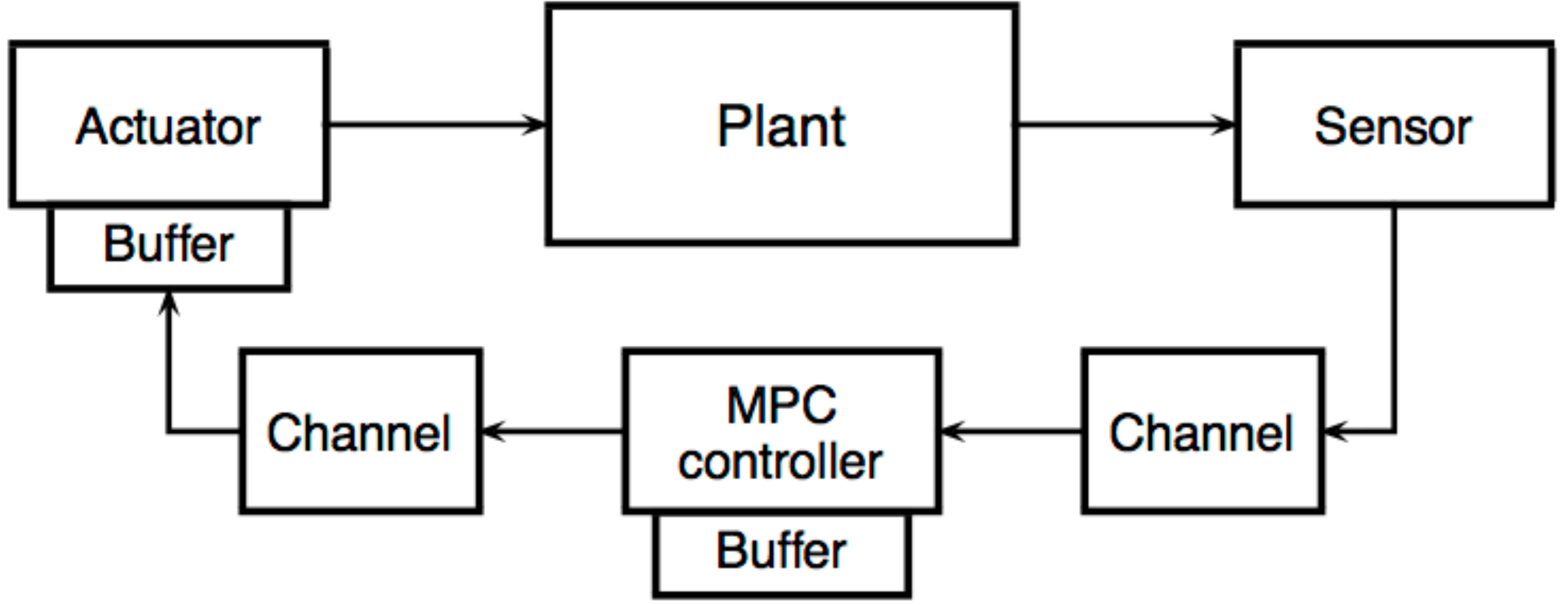}
		\caption{Scheme of the networked control system}
		\label{fig:networked system}
	\end{center}
\end{figure}
\else
\begin{figure}[!ht]
	\centering
	\psset{xunit=0.75cm,yunit=0.75cm,runit=0.75cm,arrowscale=1}
	\linethickness{0.5pt}
	\begin{pspicture}(0.0,-2.5)(10.0,1.0)
		\sffamily
		\linethickness{1pt}
		\put(0.0,-0.5){\framebox(2.0,1.0)}
		\rput[0](1.0,0.0){\parbox[]{2.0cm}{\centering \relsize{-2}{Actuator}}}
		\put(0.25,-1.05){\framebox(1.5,0.5)}
		\rput[0](1.0,-0.76){\parbox[]{2.0cm}{\centering \relsize{-2}{Buffer}}}
		\psline{->}(2.0, 0.0)(3.5, 0.0)
		\put(3.5,-0.75){\framebox(3.0,1.5)}
		\rput[0](5.0,0.0){\parbox[]{2.0cm}{\centering \relsize{-1}{Plant}}}
		\psline{->}(6.5, 0.0)(8.0, 0.0)
		\put(8.0,-0.5){\framebox(2.0,1.0)}
		\rput[0](9.0,0.0){\parbox[]{2.0cm}{\centering \relsize{-2}{Sensor}}}
		\psline{->}(9.0, -0.5)(9.0, -2.0)(8.5,-2.0)
		\put(7.0,-2.5){\framebox(1.5,1.0)}
		\rput[0](7.75,-2.0){\parbox[]{1.5cm}{\centering \relsize{-2}{Channel}}}
		\psline{->}(7.0, -2.0)(6.0, -2.0)
		\put(4.0,-2.5){\framebox(2.0,1.0)}
		\rput[0](5.0,-2.0){\parbox[]{2.0cm}{\centering \relsize{-2}{MPC\\controller}}}
		\put(4.2,-3.05){\framebox(1.6,0.5)}
		\rput[0](5.0,-2.76){\parbox[]{2.0cm}{\centering \relsize{-2}{Buffer}}}
		\psline{->}(4.0, -2.0)(3.0, -2.0)
		\put(1.5,-2.5){\framebox(1.5,1.0)}
		\rput[0](2.25,-2.0){\parbox[]{1.5cm}{\centering \relsize{-2}{Channel}}}
		\psline{->}(1.5, -2.0)(1.0, -2.0)(1.0, -1.05)
	\end{pspicture}
	\caption{Scheme of the networked control system}
	\label{fig:networked system}
\end{figure}
\fi
The idea of introducing time stamps is the following: Given the latest measurement time instant $t_s$, we can use the synchronized clocks to compute the transmission delay $\tsc$ from sensor to controller. Now, in order to be implementable, we require that $\updateinstant$ in the computation of $\mu_{T, \Delta, \sigma}$ is chosen such that the control function $\mu_{T, \Delta, \sigma}$ arrives at the actuator buffer at a time $t_a \le \sigma$. To accomplish this, we need to know the computing and transmission delays $\tsc$, $\tc$ and $\tca$. Since the latter delays are not known at computation time, bounds $\tc^{\max}$ and $\tca^{\max}$ are imposed which gives us $\updateinstant = t_s + \tsc + \tc^{\max} + \tca^{\max}$, see also Fig. \ref{fig:time} for a schematical sketch.
\ifpdf
\begin{figure}[!ht]
	\begin{center}
		\includegraphics[width=0.3\textwidth]{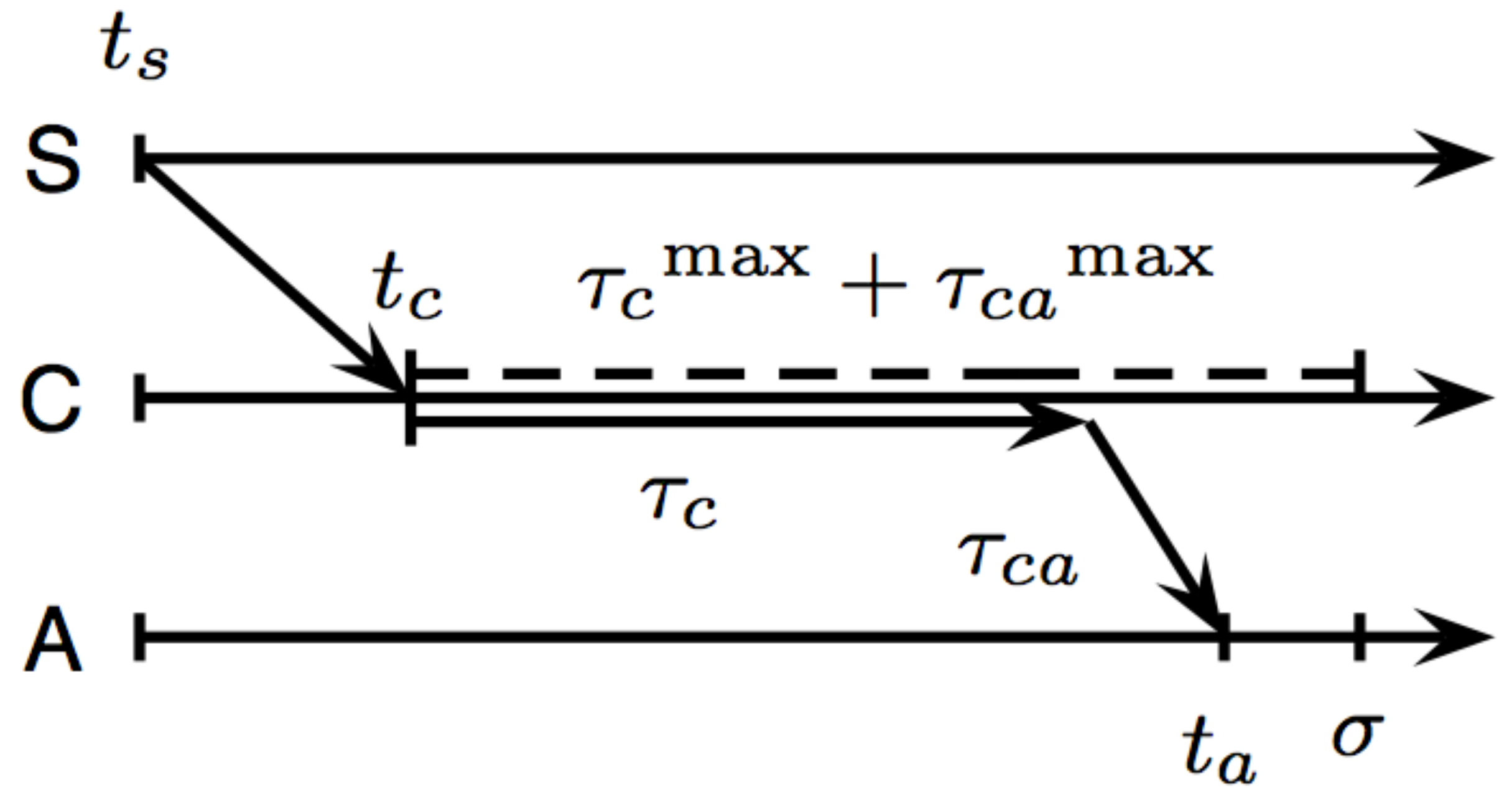}
		\caption{Timeline connections within the scheme}
		\label{fig:time}
	\end{center}
\end{figure}
\else
\begin{figure}[!ht]
\centering
    \psset{xunit=0.85cm,yunit=0.75cm,runit=0.8cm,arrowscale=1.5}
    \begin{pspicture}(-0.5,-0.2)(5.0,2.5)
        \sffamily
        \psline[linewidth=1.0pt]{->}(0.0, 0.0)(5.0, 0.0)
        \psline[linewidth=1.0pt]{-}(0.0, -0.1)(0.0, 0.1)
        \uput[180](0.0, 0.0){\relsize{-1}{A}}
        \psline[linewidth=1.0pt]{->}(0.0, 1.0)(5.0, 1.0)
        \psline[linewidth=1.0pt]{-}(0.0, 0.9)(0.0, 1.1)
        \uput[180](0.0, 1.0){\relsize{-1}{C}}
        \psline[linewidth=1.0pt]{->}(0.0, 2.0)(5.0, 2.0)
        \psline[linewidth=1.0pt]{-}(0.0, 1.9)(0.0, 2.1)
        \uput[180](0.0, 2.0){\relsize{-1}{S}}
        \uput[90](0.0, 2.1){\relsize{-1}{$t_s$}}
        \psline[linewidth=1.0pt]{->}(0.0, 2.0)(1.0, 1.0)
        \psline[linewidth=1.0pt]{-}(1.0, 0.8)(1.0, 1.2)
        \uput[90](1.0, 1.1){\relsize{-1}{$t_c$}}
        \psline[linewidth=1.0pt]{->}(1.0, 0.9)(3.5, 0.9)
        \uput[270](2.0, 0.9){\relsize{-1}{$\tc$}}
        \psline[linestyle=dashed, linewidth=1.0pt]{-}(1.0, 1.1)(3.25, 1.1)
        \uput[90](2.8, 1.1){\relsize{-1}{$\tc^{\max}+\tca^{\max}$}}
        \psline[linewidth=1.0pt]{->}(3.5, 0.9)(4.0, 0.0)
        \uput[200](3.7, 0.45){\relsize{-1}{$\tca$}}
        \psline[linestyle=dashed, linewidth=1.0pt]{-}(3.25, 1.1)(4.5, 1.1)
        \psline[linewidth=1.0pt]{-}(4.5, 1.0)(4.5, 1.2)
        \uput[270](4.0, -0.1){\relsize{-1}{$t_a$}}
        \psline[linewidth=1.0pt]{-}(4.0, -0.1)(4.0, 0.1)
        \uput[270](4.5, -0.1){\relsize{-1}{$\sigma$}}
        \psline[linewidth=1.0pt]{-}(4.5, -0.1)(4.5, 0.1)
    \end{pspicture}
\caption{Timeline connections within the scheme}
\label{fig:time}
\end{figure}
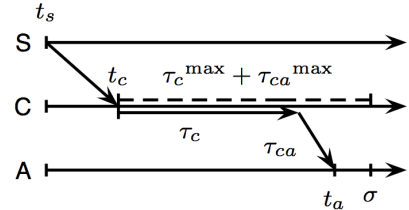
\fi
Now, we use a model based predictor within our controller as in \karl{\cite{Bemporad1998, PolushinLiuLung2008, TangSilva2007, PinParsini2008}} to compensate for these delays. If $\tc + \tca > \tc^{\max} + \tca^{\max}$ \karl{holds, the packet is considered} to be lost.

While the actuator buffer allows us to store $\mu_{T, \Delta, \sigma}$ in order to compute the control input
\begin{align}
	\label{ActuatorControl}
	u(t) = \feedback_{T, \Delta, \updateinstant}(t - \updateinstant, \state(\updateinstant))
\end{align}
the role of the controller buffer is to establish the following consistency property:
\begin{definition}
	(i) A feedback control $\feedback_{T, \Delta, \updateinstant}(\cdot, \state(\updateinstant))$ is \textit{consistently predicted} if the control $\tilde{u}_{[t_s, \updateinstant)}$ used for the prediction of $\state(\updateinstant)$ in the controller is identical to the control \juergen{$u_{[t_s, \updateinstant)}$} applied by the actuator.\\
	(ii) A networked control scheme is \textit{prediction consistent} if at each time $t \in \R$ the applied control in the actuator \eqref{ActuatorControl} is well defined, i.e. $t \leq \updateinstant + \Delta$, \karl{and }
	$\feedback_{T, \Delta, \updateinstant}$ is consistently predicted.
\end{definition}
The concept of prediction consistency allows us to separate the \textit{analysis of the feedback law}, e.g., in terms of stability and performance, from the \textit{influence of the network} on the closed--loop, e.g., in terms of robustness, cf.~\cite{FindeisenGruenePannekVarutti2010}. Examples for prediction consistent network arcitectures can be found, e.g., in \cite{Bemporad1998, PolushinLiuLung2008, GruenePannekWorthmann2009b}.

Note that MPC is ideally suited for such an application. For one, a model of the plant is already at hand, i.e. the prediction of the state measurement $\xh$ can be don\karl{e b}y evaluating the system dynamics \eqref{NotationSystemDynamics} using the control stored in the controller buffer. Secondly, the outcome of the MPC algorithm is already an optimal control defined on a control horizon of length $T$. Hence, for $\Delta\le T$ a control signal of length $\Delta$ is readily available for transmission. 

The delays $\tsc$ and $\tca$ and thus the length $\delta$ of the portion of each $\mu_{T, \Delta, \sigma}$ used at the actuator may vary with time depending on the current network load. In principle, one could make this value independent of time by setting $\delta=\Delta$, however, for robustness reasons it is desirable to always use the most recent control input, i.e., to keep $\delta$ as small as possible \cite{FindeisenGruenePannekVarutti2010}.

Hence, time varying control horizons $\delta$ should be considered and our goal is thus to find a condition which ensures asymptotic stability of the MPC closed--loop in this situation. To this end, we aim at employing the value function $V_T(\cdot)$ as a common Lyapunov function and show the relaxed Lyapunov inequality
\begin{equation}\label{NotationRelaxedLyapunovInequalityContinuousTime}
	V_T(\x_{\mu_{T}}(\delta;\xh)) \leq V_T(\xh) - \alpha \int_0^{\delta} \l(\x_{\mu_{T}}(t;\xh),\mu_{T}(t,\xh))\, dt
\end{equation}
with $\alpha \in (0,1]$ for each feasible state $\x \in \Xb$ and all $\delta \in [\tc,T-\tc]$, \karl{cf.~
\cite{GruenePannekSeehaferWorthmann2010}} for a discrete settin\karl{g.} 

\section{Recap}\label{SectionRecap}

Our main tool in order to establish the relaxed Lyapunov Inequality \eqref{NotationRelaxedLyapunovInequalityContinuousTime} is \karl{a} stability condition introduced in \cite{RebleAllgoewer2011} for fixed control horizon $\delta$. \karl{Assumption \ref{AssumptionControllabilityContinuousTime} is needed to state the respective result.}
\begin{assumption}\label{AssumptionControllabilityContinuousTime}
	Let $C \geq 1$ and $\mu > 0$ be given. Suppose that, for each $\x \in \Xb$ \karl{and $t \in \R_{\geq 0}$}, a control function $\karl{\u_{\x}} \in \Uc_{\x}[0,\infty)$ exists which satisfies
	\begin{equation}
		\l(\x_{\u_{\x}}(t;\x), \u_{\x}(t)) \leq C e^{-\mu t} \karl{\min_{\u \in \Ub} \l(\x,\u) =: C e^{-\mu t} \l^\star(\x)}. \nonumber
	\end{equation}
\end{assumption}

Assumption \ref{AssumptionControllabilityContinuousTime} is an exponential controllability condition in terms of the stage cos\karl{t 
w}ith overshoo\karl{t $C$} and decay rate $\mu$. \karl{Then} the main result deduced in \cite{RebleAllgoewer2011} reads as follows:
\begin{theorem}\label{TheoremAlphaContinuousTime}
	Suppose that \karl{$\eta(\|\x - \x^\star\|) \leq \l^\star(\x) $, $\x \in \Xb$, with $\eta \in \mathcal{K}_\infty$ and Assumption \ref{AssumptionControllabilityContinuousTime} %
	}hold. Furthermore, let $T > \delta > 0$ and $\overline{\alpha} \in (0,1)$ be chosen such that $\alpha_{T,\delta}$ given by
	\begin{equation}\label{TheoremAlphaContinuousTimeInequality}
		1 - \frac {\big(e^{\mu \delta} \hspace*{-.5mm} - \hspace*{-.5mm} 1\big)^{\frac 1C}} {\big(e^{\mu T} \hspace*{-.5mm} - \hspace*{-.5mm} 1\big)^{\frac 1C} \hspace*{-1.mm} - \hspace*{-.5mm} \big(e^{\mu \delta} \hspace*{-.5mm} - \hspace*{-.5mm} 1\big)^{\frac 1C}} \cdot \frac {\big(e^{\mu (T - \delta)} \hspace*{-.5mm} - \hspace*{-.5mm} 1\big)^{\frac 1C}} {\big(e^{\mu T} \hspace*{-.5mm} - \hspace*{-.5mm} 1\big)^{\frac 1C} \hspace*{-1.mm} - \hspace*{-.5mm} \big(e^{\mu (T - \delta)} \hspace*{-.5mm} - \hspace*{-.5mm} 1\big)^{\frac 1C}}
	\end{equation}
	satisfies the condition $\alpha_{T,\delta} \geq \overline{\alpha}$. Then, for each $\x_0 \in \Xb$, the MPC closed--loop solution $\x_{\mu_{T,\delta}}^{\text{MPC}}(\cdot;\x_0)$ is asympto\-tically stable and satisfies the suboptimality bound
	\begin{equation}\label{TheoremAlphaContinuousTimeEstimate}
		J_\infty^{\text{MPC}} \hmpc (\x_0) \hspace*{-0.25mm} = \hspace*{-1.25mm} \int_{0}^{\infty} \hspace*{-2.25mm} \l(\x_{\mu_{T,\delta}}^{\text{MPC}} \hmpc (t;\x_0),\mu_{T,\delta}^{\text{MPC}} \hmpc (t;\x_0))\hspace*{0.25mm} dt \leq \hspace*{-0.25mm} \frac {V_{\hinfty \infty}(\x_0)}{\overline{\alpha}} \hspace*{-0.25mm}.
	\end{equation}
\end{theorem}

Inequality \eqref{TheoremAlphaContinuousTimeEstimate} gives a performance estimates which compares the resulting MPC closed--loop costs with the theoretically achievable minimal costs on the infinite time horizon. Here, the monotonicity of the optimal value function $V_T(\cdot)$ in the prediction horizon $T$ --- an inherent property of unconstrained MPC schemes --- is crucial in order to deduce this bound on the, in general, unknown quantity $V_\infty(\cdot)$. Note that, for given suboptimality index $\overline{\alpha} \in (0,1)$, the stability condition $\alpha_{T,\delta} \geq \overline{\alpha}$ always holds for a sufficiently large prediction horizon $T$, cf.~\cite[Section 4.1]{RebleAllgoewer2011}\karl{.}


\section{Results}\label{SectionResults}

As pointed out in the previous section, in the networked context we would like to have a stability criterion for time varying control horizon $\delta$. Unfortunately, the stability condition presented in the previous Section \ref{SectionRecap} assumes a fixed control horizon. In this section we show how to extend the result to varying $\delta$. In particular, we show that the inequality $\alpha_{T,\delta} \geq \overline{\alpha}$ for some control horizon $\delta \in (0,T)$ implies $\alpha_{T,\delta} \geq \overline{\alpha}$ for every $\delta \in [\min \{\delta, T-\delta\},\max \{\delta,T-\delta\}]$. This result is the key to derive a stability theorem for time varying control horizons without imposing additional assumptions.

To this end, we show symmetry and monotonicity properties of the performance bound $\alpha_{T,\delta}$ 
with respect to the control horizon $\delta$. Indeed, symmetry is a direct consequence of the presented formula.
\begin{corollary}\label{CorollarySymmetry}
	 The performance \karl{index} $\alpha_{T,\delta}$ given by Formula \eqref{TheoremAlphaContinuousTimeInequality} satisfies $\alpha_{T,\delta} = \alpha_{T,T-\delta}$, i.e. $\alpha_{T,\delta}$ is symmetric with symmetry axis $\delta = T/2$.
\end{corollary}

In contrast to the symmetry property given in Corollary \ref{CorollarySymmetry}, deducing monotonicity of Formula \eqref{TheoremAlphaContinuousTimeInequality} in $\delta$ on $(0,T/2]$ is not straightforward. To this end, we apply results from \cite{GruenePannekSeehaferWorthmann2010} and \cite{WorthmannRebleGrueneAllgoewer2012}. In the first of these references a discrete time counterpart of Theorem \ref{TheoremAlphaContinuousTime} was presented based on a corresponding version of Assumption \ref{AssumptionControllabilityContinuousTime}. In the second article, the relation between both the discrete and continuous time approaches was investigated. To this end, so called iterative refinements were used, i.e. a sequence of discretizations such that \karl{each element }
is a partition of its predecessor. Then, monotone convergence of the corresponding suboptimality bounds to $\alpha_{T,\delta}$ given by Formula \eqref{TheoremAlphaContinuousTimeInequality} was shown for the discretization parameter tending to zero. We combine these results to show that monotonicity carries over from the discrete to the continuous time setting.
\begin{proposition}\label{PropositionMonotonicity}
	For $\delta \in (0,T/2)$, the performance estimate $\alpha_{T,\delta}$ given by Formula \eqref{TheoremAlphaContinuousTimeInequality} has the monotonicity property
	\begin{equation}\label{PropositionMonotonicityInequality}
		\alpha_{T,\delta} \leq \alpha_{T,\tilde{\delta}} \qquad\text{for }\tilde{\delta} \in [\delta,T/2].
	\end{equation}
\end{proposition}\hv
\proof 
We show \eqref{PropositionMonotonicityInequality} for $T, \z, \zt \in\Q$. The assertion for real values $T, \z, \zt \in\R$ then follows from the fact that $\Q$ is dense in $\R$ and that $\alpha_{T,\delta}$ is continuous in $T$ and $\delta$.

	Consider $T, \z, \zt\in \Q$ with $\zt \in [\z,T/2]$. Since $\z$ and $\zt$ are rational numbers, there exist $\s \in \Q$ and $\m, \mt, N \in \N$ such that $\m \s = \z$, $\mt \s = \zt$ $N \s = T$ hold, i.e. $\s$ is a common denominator. Then, using the abbreviations $N_k := 2^k N$, $\m_k := 2^k \m$, we define $\alpha_k = \alpha_k(N,\m)$ by
	\begin{equation}\label{PropositionMonotonicityProofAlphaFormulaDiscrete}
		1 - \frac {\prod\limits_{m_k+1}^{N_k} \hspace*{-1.mm} (\gamma_i^k - 1) \prod\limits_{N_k-m_k+1}^{N_k} \hspace*{-4mm}(\gamma_i^k - 1)}{\bigg[\prod\limits_{m_k+1}^{N_k} \hspace*{-1.25mm} \gamma_i^k - \hspace*{-2.25mm} \prod\limits_{m_k+1}^{N_k} \hspace*{-1.mm} (\gamma_i^k - 1) \hspace*{-0.5mm} \bigg]  \hspace*{-1.mm} \bigg[\prod\limits_{N_k-m_k+1}^{N_k} \hspace*{-4.5mm} \gamma_i^k \hspace*{2mm} - \hspace*{-2mm} \prod\limits_{N_k-m_k+1}^{N_k} \hspace*{-4.mm} (\gamma_i^k - 1) \hspace*{-0.5mm} \bigg]}
	\end{equation}
	with $\gamma_i^k := C \sum_{n=0}^{i-1} \sigma^{2^{-k}}$ with $\sigma := e^{- \mu \s} \in (0,1)$. Analogously, $\tilde{\alpha}_k = \alpha_k(N,\mt)$ is defined with $\mt_k := 2^k \mt$ instead of $\m_k$. Summarizing, $\alpha_k$ and $\tilde{\alpha}_k$ only deviate in the parameters $\m$ and $\mt$. Furthermore, note that $\m \leq \mt$ holds.

	Formula \eqref{PropositionMonotonicityProofAlphaFormulaDiscrete} represents a discrete time counterpart of Formula \eqref{TheoremAlphaContinuousTimeInequality}, cf.~\cite{WorthmannRebleGrueneAllgoewer2012} for details. In particular, $\alpha_k$ and $\tilde{\alpha}_k$ are monotonically increasing in $k$, cf.~\cite[Proposition 3.3]{WorthmannRebleGrueneAllgoewer2012}. Additionally $\alpha_k$ and $\tilde{\alpha}_k$ converge to $\alpha_{T,\z}$ and $\alpha_{T,\zt}$, respectively, cf.~\cite[Theorem 3.2]{WorthmannRebleGrueneAllgoewer2012}. Next, we can use the fact that the discrete time suboptimality bounds $\alpha_k(N,\cdot)$, $k \in \N_0$, are monotonically increasing in the parameter $\m$ and, thus, $\alpha_k = \alpha_k(N,\m) \leq \alpha_k(N,\mt) = \tilde{\alpha}_k$ holds for each $k \in \N_0$, cf.~\cite[Proposition 7.3]{GruenePannekSeehaferWorthmann2010}. Hence, this inequality also holds for the limit, i.e.
	\begin{equation*}
		\alpha_{T,\z} = \lim_{k \rightarrow \infty} \alpha_k(N,\m) \leq \lim_{k \rightarrow \infty} \alpha_k(N,\mt) = \alpha_{T,\zt}
	\end{equation*}
	which shows \eqref{PropositionMonotonicityInequality} for $T, \z, \zt \in\Q$.
\endproof
Note that $\alpha_{T,\delta}$ defined by Formula \eqref{TheoremAlphaContinuousTimeInequality} is differentiable on $(0,T)$. From Proposition \ref{PropositionMonotonicity} we can thus conclude $\frac {\partial \alpha_{T,\delta}}{\partial \delta} \geq 0$ for all $\delta \in (0,T/2]$ and $\frac {\partial \alpha_{T,\delta}}{\partial \delta} \le 0$ on $[T/2,T)$.
\begin{remark}
	Note that the discretization procedure used in the proof of Proposition \ref{PropositionMonotonicity} was employed in a purely theoretical fashion. Checking the assumptions on the (discrete time) dynamics from \cite{WorthmannRebleGrueneAllgoewer2012} is not needed since we only use properties of Formula \eqref{PropositionMonotonicityProofAlphaFormulaDiscrete} which are independent of its connection to discrete time dynamics.
\end{remark}

Combining Corollary \ref{CorollarySymmetry} and Proposition \ref{PropositionMonotonicity} yields the following theorem as a direct consequence.
\begin{theorem}
	Let Assumption \ref{AssumptionControllabilityContinuousTime} and $\alpha_{T,\underline{\delta}} \geq \overline{\alpha} > 0$ hold for a minimal control horizon $\underline{\delta} \in (0,T)$. Then, $\alpha_{T,\delta} \geq \overline{\alpha}$ holds for all $\delta \in [\underline{\delta},T-\underline{\delta}]$. Hence, if the stability condition $\alpha_{T,\underline{\delta}} \geq \overline{\alpha}$ holds for a desired guaranteed performance bound $\overline{\alpha}$, then at least the same performance can be guaranteed for all control horizons $\delta \in [\underline{\delta},T-\underline{\delta}]$ and, thus, also for time varying control horizons $(\delta_i)_{i \in \N_0} \subset [\underline{\delta},T-\underline{\delta}]$.
\end{theorem}
\proof
	Combining the properties shown in Corollary \ref{CorollarySymmetry} and Proposition \ref{PropositionMonotonicity} immediately implies the claimed inequality $\alpha_{T,\delta} \geq \alpha_{T,\underline{\delta}} \geq \overline{\alpha}$. The assertion for time varying control horizons can be concluded analogously to the stability theorem for fixed \karl{$\delta$}, cf.~\cite{RebleAllgoewer2011}, with the same modifications being carried out in \cite{GruenePannekSeehaferWorthmann2010} to extend the discrete time stability theore\karl{m 
	f}rom fixed to time varying control horizon \karl{($V_T(\cdot)$ is used as a common Lyapunov function)}.
\endproof
\begin{remark}
	It is possible to replace the exponential controllability Assumption \ref{AssumptionControllabilityContinuousTime} by alternative controllability conditions. In this case, however, in the continuous time setting no closed formulas for $\alpha_{T,\delta}$ as in Theorem \ref{TheoremAlphaContinuousTime} are known and thus checking monotonocity may become a difficult task.
\end{remark}
\begin{remark}[Sampled-data systems]
	If a sampled-data implementation of system \eqref{NotationSystemDynamics} is necessary, the methodology proposed in \cite{WorthmannRebleGrueneAllgoewer2012} can be used to determine the required sampling rate such that the continuous time estimate is approximated arbitrarily well. This allows to transfer our results to a sampled-data framewor\karl{k.} 
\end{remark}

So far, we showed that using larger control horizons is not harmful with respect to the performance. Next, for a given prediction horizon, we compute parameter combinations $(C,\mu)$ for which stability is ensured by Theorem \ref{TheoremAlphaContinuousTime} in dependence of the control horizon $\delta$, cf.~Fig.~\ref{FigureStabilityRegionControlHorizon}.
\begin{figure}[thpb]
	\begin{center}
	  \includegraphics[width=6cm,height=5cm]{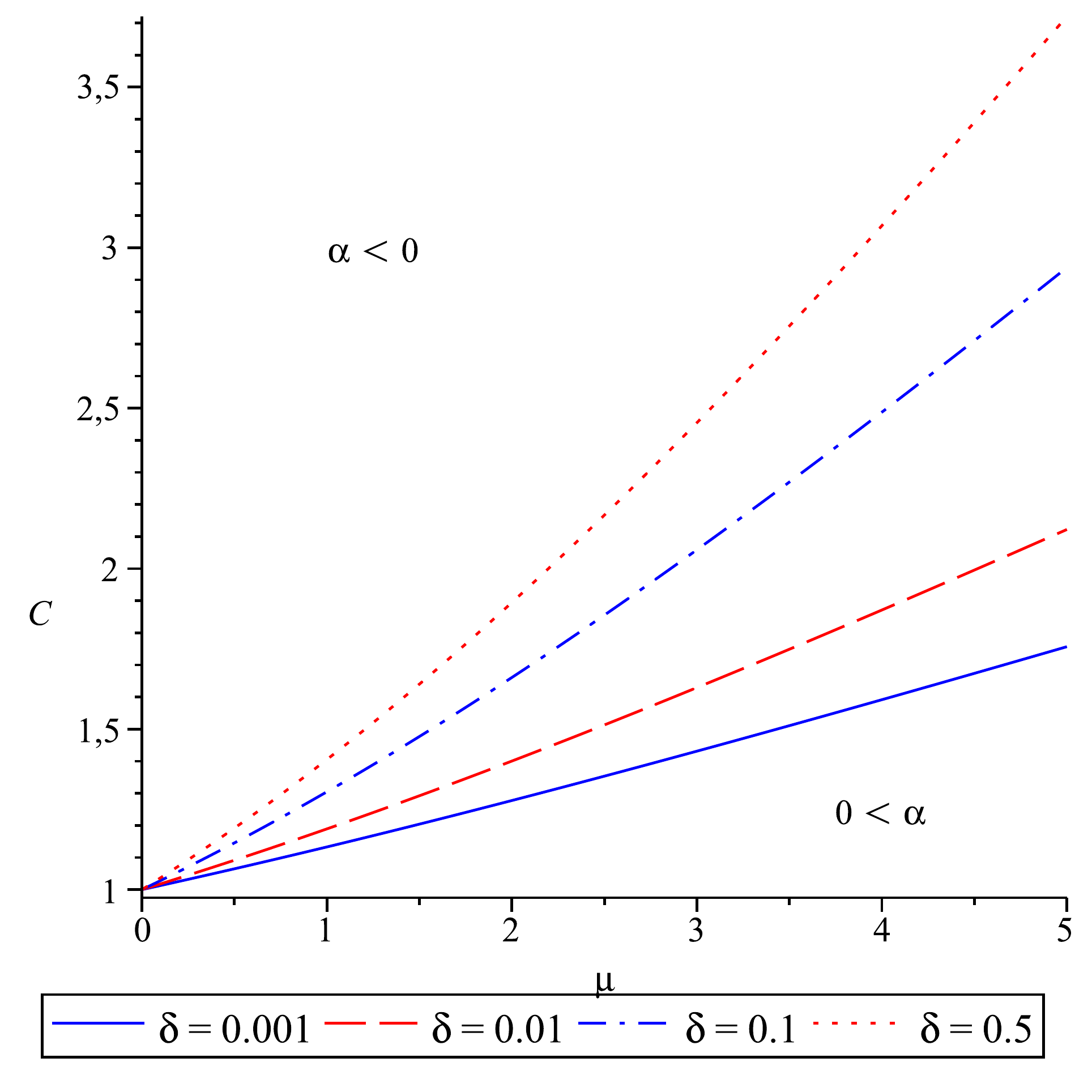}
	  \caption{Parameter pairs $(C,\sigma)$ for which $\alpha_{T,\delta} \geq 0$ is ensured by Theorem \ref{TheoremAlphaContinuousTime} for $T = 1$.}
		\label{FigureStabilityRegionControlHorizon}
	\end{center}
\end{figure}

\juergen{The respective} structurally different impact of the overshoot bound $C$ and the decay rate $\mu$ from Assumption \ref{AssumptionControllabilityContinuousTime} is investigated \juergen{in the next two propositions}.
\begin{proposition}\label{PropositionOvershoot}
	Let a decay rate $\mu > 0$ and times $T\ge \delta>0$ be given. Then, if Assumption \ref{AssumptionControllabilityContinuousTime} is satisfied with a sufficiently small overshoot bound $C \geq 1$ the stability condition $\alpha_{T,\delta} > 0$ is satisfied and, as a consequence, asymptotic stability of the MPC closed--loop is guaranteed.
\end{proposition}
\proof
	Let $C$ be equal to $1$. Then, $\alpha_{T,\delta}$ from \eqref{TheoremAlphaContinuousTimeInequality} equals
	\begin{eqnarray*}
		\alpha_{T,\delta} & = & 1 - \frac {(e^{\mu \delta} - 1)(e^{\mu(T-\delta)}-1)}{\left[e^{\mu\delta}(e^{\mu(T-\delta)}-1)\right]\left[e^{\mu(T-\delta)}(e^{\mu\delta}-1)\right]} \\
		& = & 1 - e^{-\mu T} > 0.
	\end{eqnarray*}
	Since \eqref{TheoremAlphaContinuousTimeInequality} is continuous with respect to $C$, choosing the overshoot $C > 1$ sufficiently close to one ensures $\alpha_{T,\delta} > 0$ and, thus, the assertion.
\endproof
Note that the performance estimate is bounded by $1 - e^{-\mu T}$ independent of the overshoot $C$. Hence, a sufficiently large prediction horizon $T$ may be needed in order to guarantee a desired suboptimality bound $\overline{\alpha} \in (0,1)$. \juergen{Moreover,} 
the obtained bound does not depend on the control horizon $\delta$.

In the following proposition the roles of $C$ and $\mu$ are reversed. In contrast to the previous observation that, for given decay rate $\mu$ and prediction horizon $T$, the maximal achievable performance was bounded, here, an arbitrary performance $\overline{\alpha} \in (0,1)$ can be obtained by choosing the decay rate appropriately without prolonging the prediction horizon $T$ or changing the overshoot $C$.
\begin{proposition}\label{PropositionDecayRate}
	Let an overshoot $C \geq 1$, times $T\ge \delta>0$ and a desired performance bound $\overline{\alpha} \in (0,1)$ be given. Then, if Assumption \ref{AssumptionControllabilityContinuousTime} is satisfied with a sufficiently large decay $\mu > 0$ the condition $\alpha_{T,\delta} \geq \overline{\alpha}$ is satisfied.
\end{proposition}
\proof
	Let the function $f: \R_{>0} \rightarrow \R$ be defined by $f(x) = x^{1/C}$. Then, its first derivative is given by the monotonically decreasing function
	\begin{equation*}
		f^\prime(x) = \frac {1}{C x^{1-1/C}}.
	\end{equation*}
	Consequently, for real numbers $\bar{x} \geq x > 1$, the inequality $f(x) - f(x-1) \geq f(\bar{x}) - f(\bar{x}-1)$ holds which implies
	\begin{eqnarray*}
		(e^{\mu T} \hspace*{-1.mm} - \hspace*{-.75mm} 1)^{1/C} \hspace*{-1.5mm} - (e^{\mu\delta} \hspace*{-1.mm} - \hspace*{-.75mm} 1)^{1/C} & \hspace*{-2.mm} \geq & \hspace*{-2.25mm} (e^{\mu T})^{1/C} \hspace*{-1.5mm} - (e^{\mu\delta})^{1/C}\hspace*{-1.mm}, \\
		(e^{\mu T} \hspace*{-1.mm} - \hspace*{-.75mm} 1)^{1/C} \hspace*{-1.5mm} - (e^{\mu(T-\delta)} \hspace*{-1.mm} - \hspace*{-.75mm} 1)^{1/C} & \hspace*{-2.mm} \geq & \hspace*{-2.25mm} (e^{\mu T})^{1/C} \hspace*{-1.5mm} - (e^{\mu(T-\delta)})^{1/C}\hspace*{-1.mm}.
	\end{eqnarray*}
	Using these inequalities and the fact that the nominator of the second summand in Formula \eqref{TheoremAlphaContinuousTimeInequality} is smaller than $(e^{\mu\delta})^{1/C} (e^{\mu(T-\delta)})^{1/C} = (e^{\mu T})^{1/C}$ leads to the estimate
	\begin{eqnarray*}
		\alpha_{T,\delta} & \hspace*{-1mm} \geq & \hspace*{-0.5mm} 1 - \frac {(e^{\mu T})^{1/C}}{(e^{\mu T})^{1/C} \left[(e^{\mu(T-\delta)})^{1/C}-1\right] \left[(e^{\mu \delta})^{1/C}-1\right]} \\
		& \hspace*{-1mm} = & \hspace*{-0.5mm} 1 - \left( \left[(e^{\mu(T-\delta)})^{1/C}-1\right] \left[(e^{\mu \delta})^{1/C}-1\right] \right)^{-1}.
	\end{eqnarray*}
	Since the right hand side of this expression converges to one for $\mu$ approaching infinity, a (sufficiently large) decay rate $\mu$ exists such that $\alpha_{T,\delta} \geq \overline{\alpha}$ is ensured.
\endproof
Propositions \ref{PropositionOvershoot} and \ref{PropositionDecayRate} provide some insight into the structure of the proposed controllability Assumption \ref{AssumptionControllabilityContinuousTime}. The integral $\int_0^T C e^{- \mu t}\, dt$ is the destinctive feature. This quantity converges to zero for $\mu$ approaching infinity. For given $\mu$, however, $\int_0^T e^{- \mu t}\, dt$ is a lower bound. This explains why, independently of the overshoot $C$, the performance guarantee cannot be arbitrarily well. Note that Proposition \ref{PropositionDecayRate} cannot be obtained in a discrete time setting \cite{GruenePannekSeehaferWorthmann2010}.

\section{Example}\label{SectionExample}

We illustrate our results by computing the $\alpha$-values in the relaxed Lyapunov Inequality \eqref{NotationRelaxedLyapunovInequalityContinuousTime} along simulated trajectories \juergen{for various control horizon sequences and comparing} them with our theoretical findings.
We consider a continuously stirred tank reactor (CSTR) with energy balance and reaction $A \to B$ given by the dynamics
\begin{align*}
	\dot{x}_1 & = \frac{q (x_1^f - x_1)}{V} - k_0 x_1 e^{-E_a / x_2} \\
	\dot{x}_2 & = \frac{q (x_2^f - x_2)}{V} + \frac{h}{\rho \cdot c} k_0 x_1 e^{-E_a / x_2} + \frac{\alpha}{V \cdot \rho \cdot c} (u - x_2)
\end{align*}
where the states $x_1$ and $x_2$ denote the concentration in $\frac{\text{mol}}{\text{m}^3}$ and the temperature in $\text{K}$, and the control $u$ corresponds to the temperature of the cooling jacket in $\text{K}$, cf.~\cite{HensonSeborg1997, HicksRay1971}. 
\begin{table}[!ht]
	\begin{center}
		\begin{tabular}{|c|c|c|c|} \hline
			name & symbol & quantity & unit\\ \hline\hline
			flowrate & $q$ & $100$ & $m^3/sec$ \\
			CSTR volume & $V$ & $100$ & $m^3$ \\
			pre-exponential factor & $k_0$ & $7.2\cdot10^{10}$ & $1/sec$\\
			activ.\ energy [gas const.] & $E_a$ & $8750$ & $1/K$\\
			heat of the reaction & $h$ & $5\cdot10^4$ & $J/mol$\\
			$A/B$ mixture density & $\rho$ & $1000$ & $kg/m^3$\\
			$A/B$ mixture heat capacity & $c$ & $0.239$ & $J/(kg\cdot K)$\\
			heat transfer & $\alpha$ & $5\cdot10^4$ & $W/K$\\
			feed concentration & $x_1^f$ & $1$ & $mol/m^3$\\
			feed temperature & $x_2^f$ & $350$ & $K$\\
		\hline
		\end{tabular}
		\caption{Parameters for the CSTR example.}
		\label{TableCSTR}
	\end{center}
\end{table}

Here, \karl{the} aim is to stabilize the equilibrium $\statestar = ( 0.5, 350 )$ and $\controlstar = 300$. The system is subject to \karl{the physically motivated} state constraints $\constraintstateset = [0, 1] \times [0, \infty)$. \kawo{Additionally, control} constraints $\constraintcontrolset = [250, 450]$ are imposed. The control $\control$ is considered to be piecewise constant on intervals of length $\delta = 10^{-2}$. The optimal control in each MPC step is computed over an optimization horizon $T = 0.3$ (although the CSTR can be practically stabilized using MPC with $T = 0.02$, we use the larger $T$ in our simulations in order to have more flexibility in choosing \karl{$\delta$}).

We use the initial value $\state_0 = ( 0.35, 370 )$ and define the running costs via
\begin{align*}
	\stagecost(\state, \control) = & \left( \frac{\state^\star_2}{\state^\star_1} \right)^2 ( \state - \state^\star_1 )^2 + ( \state - \state^\star_2)^2 + 10^{-3} ( \control - \controlstar)^2
\end{align*}
which \kawo{render} the state components equally weighted. The error tolerance for the used optimization and integration method were set to $10^{-7}$ and \kawo{$10^{-6}$, respectively}. 
Moreover, we added a truncation region of the stage costs $\stagecost$ of $\varepsilon = 10^{-12}$ to compensate for a possible practical stability region and numerical errors, cf. \cite[Theorem 21]{GruenePannek2009} for details. The performance index $\alpha_{T, \delta}$ can then be evaluated along the simulated closed--loop trajectory via
\begin{equation}\label{NotationSuboptimalityEstimateAlongTrajectory}
	\alpha_{T, \delta} = \inf_{n \in \{n | \exists k \in \N_0: n = k \delta \}} \alpha_{T, \delta}(n)
\end{equation}
where the local performance index $\alpha_{T, \delta}(n)$ is given by
\begin{align*}
	\alpha_{T, \delta}(n) = \frac{V_T(\x_{\mu}(n \delta; \x_0)) - V_T(\x_{\mu}((n+1) \delta; \x_0))}{\int\limits_{0}^{\delta} (\stagecost(\x_{\mu}(t; \x_{\mu}(n \delta; \x_0)),\mu(t;\x_{\mu}(n \delta; \x_0)))) dt  - \varepsilon }
\end{align*}
with $\mu = \mu_{T, \delta}$ if the denominator of the right hand side is strictly positive and $\alpha_{T, \delta}(n) = 1$ otherwise. Note that $\alpha_{T, \delta}$ \juergen{is} 
negative if the value function increases. 

\karl{Fig. \ref{fig:trajectories} shows 
the closed--loop solution 
in dependency of different \juergen{fixed and time varying control horizon sequences}}. 
Although the trajectories corresponding to $\delta \in \{0.15, 0.25,0.30\}$ appear to be less efficient due to avoiding the turning point close to $(0.4, 325)$, 
a} computation of the performance index $\alpha_{0.3, \delta}$ reveals the opposite, cf.~Fig.~\ref{fig:alpha}: If the control horizon is increased, the performance index rises as expected from our theoretical results in Section \ref{SectionResults}.
\begin{figure}[!ht]
	\begin{center}
		\includegraphics[width=0.48\textwidth]{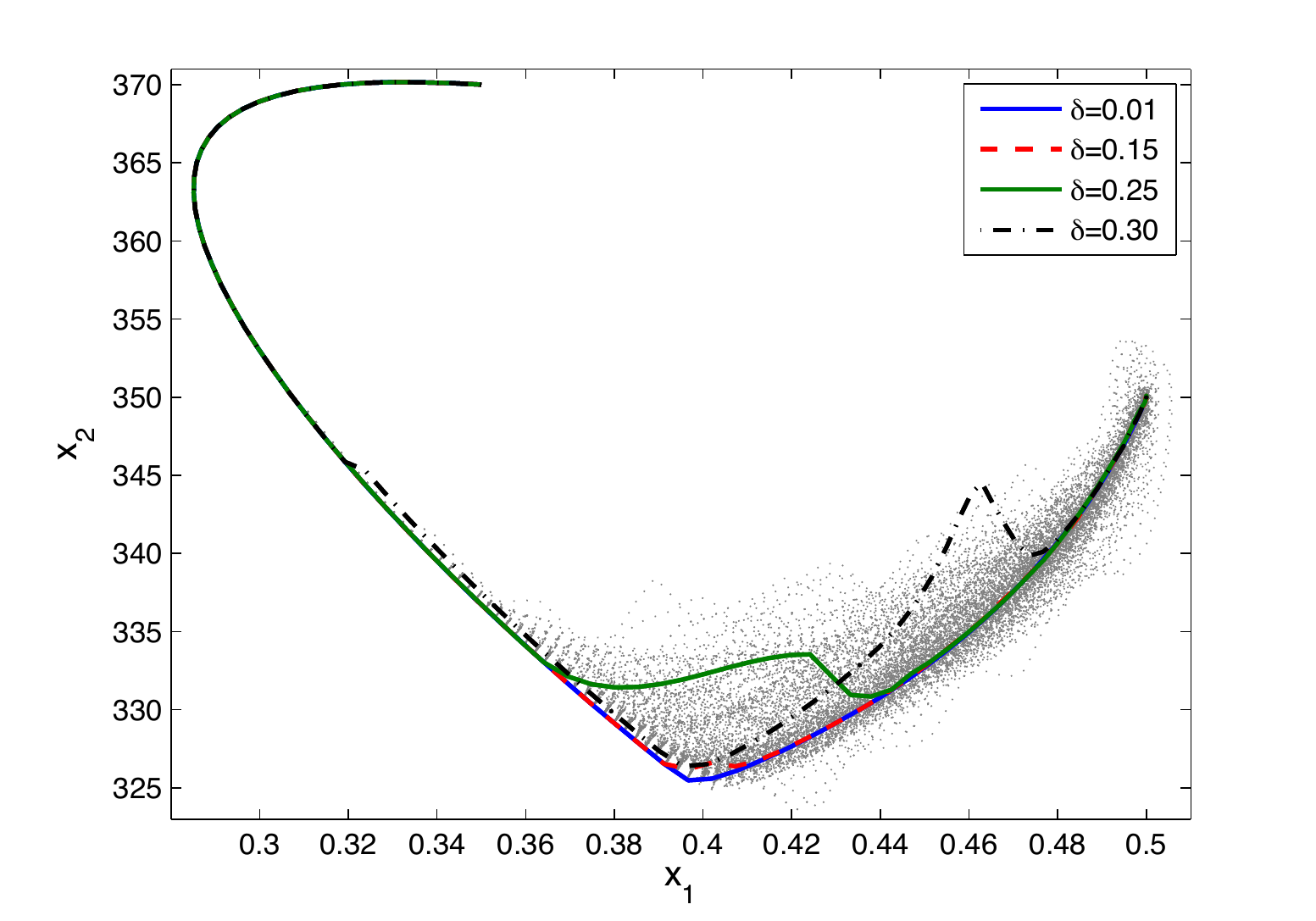}
		\caption{Closed--loop solutions for \kawo{different fixed control horizons $\delta$ and $400$ closed-loop solutions for randomly time varying $\delta \in [0.1,0.3]$.}}
		\label{fig:trajectories}
	\end{center}
\end{figure}
\begin{figure}[!t]
	\begin{center}
		\includegraphics[width=0.48\textwidth]{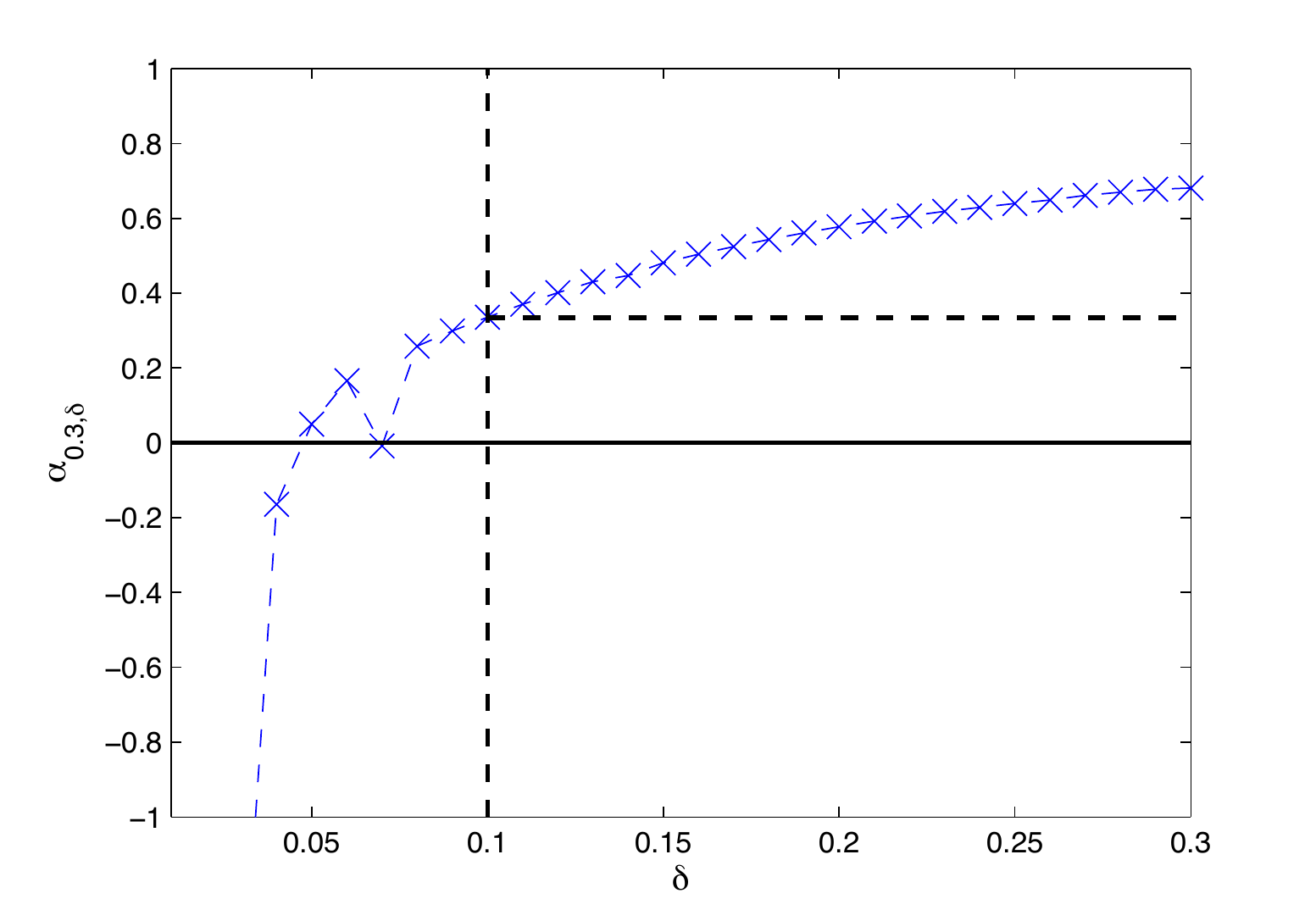}
		\caption{Development of $\alpha_{0.3, \delta}$ for \kawo{fixed} control horizons $\delta$. \kawo{$\alpha_{0.3,0.1}$ yields a lower bound for time varying control horizon $\delta \in [0.1,0.3]$.}}
		\label{fig:alpha}
	\end{center}
\end{figure}

For this particular example, \karl{symmetry of $\alpha_{T, \delta}$ with respect to the control horizon $\delta$ is not observed}. Instead, the performance index increases (almost) monotonically. There are, however, examples in which the suboptimality degree decays rapidly for longer control horizons, \karl{cf.~\cite{GruenePannekWorthmann2009ECC}}. Utilizing the outcome of Fig. \ref{fig:alpha}, we restricted $\delta$ to $ [0.1, 0.3]$ and computed the closed loop solutions shown in Fig. \ref{fig:trajectories}. From this, we \kawo{observed} that the suboptimality bound \kawo{$\alpha_{0.3, 0.1} \approx 0.3346$ (obtained for fixed control horizon $\delta = 0.1$) holds for control horizons randomly varying in a rather large interval}. Fo\kawo{r s}hort control horizons on the other hand, the relaxed Lyapunov Inequality \eqref{NotationRelaxedLyapunovInequalityContinuousTime} \karl{---  our main tool in order to ensure stability ---} is violated\karl{, cf.~Fig. \ref{fig:alpha}}. \kawo{In an a posteriori analysis, we found that the closed--loop costs $J_\infty^{\text{MPC}}(\x_0)$ are almost constant with respect to $\delta$.} In conclusion, \juergen{this example demonstrates that packet dropouts and non-negligible delays} \kawo{can} be compensated \juergen{via} \karl{time varying} control horizons without affecting our (nominal) stability and performance estimates.

\section{Conclusions}
We presented stability and performance estimates for nonlinear model predictive control with time varying control horizon in a continuous time setting. In particular, we deduced symmetry and monotonicity properties for the stability condition introduced in \cite{RebleAllgoewer2011} which were exploited in order to show that no additional assumptions are needed compared to the analysis of schemes with fixed (short) control horizons. The results can be used in order to obtain (nominal) stability and performance estimates for prediction consistent networked control schemes.

\addtolength{\textheight}{-10.3cm}   




\section*{Acknowledgement}

This work was supported by the DFG priority research program 1305 ``Control Theory of Digitally Networked Dynamical Systems''.



\end{document}